\documentclass[12pt,reqno]{amsart}
\usepackage{amscd,amssymb,amsmath,multicol,float,scalefnt,cancel}
\restylefloat{table}
\newtheorem{thm}[equation]{Theorem}
\numberwithin{equation}{section}

\newtheorem{expl}[equation]{Example}

\begin{document}
\raggedbottom \voffset=-.7truein \hoffset=0truein \vsize=8truein
\hsize=6truein \textheight=8truein \textwidth=6truein
\baselineskip=18truept

\def\mapright#1{\ \smash{\mathop{\longrightarrow}\limits^{#1}}\ }
\def\mapleft#1{\smash{\mathop{\longleftarrow}\limits^{#1}}}
\def\mapup#1{\Big\uparrow\rlap{$\vcenter {\hbox {$#1$}}$}}
\def\mapdown#1{\Big\downarrow\rlap{$\vcenter {\hbox {$\ssize{#1}$}}$}}
\def\mapne#1{\nearrow\rlap{$\vcenter {\hbox {$#1$}}$}}
\def\mapse#1{\searrow\rlap{$\vcenter {\hbox {$\ssize{#1}$}}$}}
\def\mapr#1{\smash{\mathop{\rightarrow}\limits^{#1}}}
\def\ss{\smallskip}
\def\vp{v_1^{-1}\pi}
\def\at{{\widetilde\alpha}}
\def\sm{\wedge}
\def\la{\langle}
\def\ra{\rangle}
\def\on{\operatorname}
\def\ol#1{\overline{#1}{}}
\def\spin{\on{Spin}}
\def\lbar{ \ell}
\def\qed{\quad\rule{8pt}{8pt}\bigskip}
\def\ssize{\scriptstyle}
\def\a{\alpha}
\def\bz{{\Bbb Z}}
\def\im{\on{im}}
\def\ct{\widetilde{C}}
\def\ext{\on{Ext}}
\def\sq{\on{Sq}}
\def\eps{\epsilon}
\def\ar#1{\stackrel {#1}{\rightarrow}}
\def\br{{\bold R}}
\def\bC{{\bold C}}
\def\bA{{\bold A}}
\def\bB{{\bold B}}
\def\bD{{\bold D}}
\def\bh{{\bold H}}
\def\bQ{{\bold Q}}
\def\bP{{\bold P}}
\def\bx{{\bold x}}
\def\bo{{\bold{bo}}}
\def\si{\sigma}
\def\Ebar{{\overline E}}
\def\dbar{{\overline d}}
\def\Sum{\sum}
\def\tfrac{\textstyle\frac}
\def\tb{\textstyle\binom}
\def\Si{\Sigma}
\def\w{\wedge}
\def\equ{\begin{equation}}
\def\b{\beta}
\def\G{\Gamma}
\def\g{\gamma}
\def\k{\kappa}
\def\psit{\widetilde{\Psi}}
\def\tht{\widetilde{\Theta}}
\def\psiu{{\underline{\Psi}}}
\def\thu{{\underline{\Theta}}}
\def\aee{A_{\text{ee}}}
\def\aeo{A_{\text{eo}}}
\def\aoo{A_{\text{oo}}}
\def\aoe{A_{\text{oe}}}
\def\vbar{{\overline v}}
\def\endeq{\end{equation}}
\def\sn{S^{2n+1}}
\def\zp{\bold Z_p}
\def\A{{\cal A}}
\def\P{{\mathcal P}}
\def\cj{{\cal J}}
\def\zt{{\bold Z}_2}
\def\bs{{\bold s}}
\def\bof{{\bold f}}
\def\bq{{\bold Q}}
\def\be{{\bold e}}
\def\Hom{\on{Hom}}
\def\ker{\on{ker}}
\def\coker{\on{coker}}
\def\da{\downarrow}
\def\colim{\operatornamewithlimits{colim}}
\def\zphat{\bz_2^\wedge}
\def\io{\iota}
\def\Om{\Omega}
\def\Prod{\prod}
\def\e{{\cal E}}
\def\exp{\on{exp}}
\def\abar{{\overline a}}
\def\xbar{{\overline x}}
\def\ybar{{\overline y}}
\def\zbar{{\overline z}}
\def\Rbar{{\overline R}{}}
\def\nbar{{\overline n}}
\def\cbar{{\overline c}}
\def\qbar{{\overline q}}
\def\bbar{{\overline b}}
\def\et{{\widetilde E}}
\def\ni{\noindent}
\def\coef{\on{coef}}
\def\den{\on{den}}
\def\lcm{\on{l.c.m.}}
\def\vi{v_1^{-1}}
\def\ot{\otimes}
\def\psibar{{\overline\psi}}
\def\mhat{{\hat m}}
\def\exc{\on{exc}}
\def\ms{\medskip}
\def\ehat{{\hat e}}
\def\etao{{\eta_{\text{od}}}}
\def\etae{{\eta_{\text{ev}}}}
\def\dirlim{\operatornamewithlimits{dirlim}}
\def\gt{\widetilde{L}}
\def\lt{\widetilde{\lambda}}
\def\st{\widetilde{s}}
\def\ft{\widetilde{f}}
\def\sgd{\on{sgd}}
\def\lfl{\lfloor}
\def\rfl{\rfloor}
\def\ord{\on{ord}}
\def\gd{{\on{gd}}}
\def\rk{{{\on{rk}}_2}}
\def\nbar{{\overline{n}}}
\def\lg{{\on{lg}}}
\def\cR{\mathcal{R}}
\def\cS{\mathcal{S}}
\def\cT{\mathcal{T}}
\def\N{{\Bbb N}}
\def\Z{{\Bbb Z}}
\def\Q{{\Bbb Q}}
\def\R{{\Bbb R}}
\def\C{{\Bbb C}}
\def\l{\left}
\def\r{\right}
\def\mo{\on{mod}}
\def\xt{\times}
\def\notimm{\not\subseteq}
\def\Remark{\noindent{\it  Remark}}

\def\*#1{\mathbf{#1}}
\def\0{$\*0$}
\def\1{$\*1$}
\def\22{$(\*2,\*2)$}
\def\33{$(\*3,\*3)$}
\def\ss{\smallskip}
\def\ssum{\sum\limits}
\def\dsum{\displaystyle\sum}
\def\la{\langle}
\def\ra{\rangle}
\def\on{\operatorname}
\def\o{\on{o}}
\def\U{\on{U}}
\def\lg{\on{lg}}
\def\a{\alpha}
\def\bz{{\Bbb Z}}
\def\eps{\varepsilon}
\def\bc{{\bold C}}
\def\bN{{\bold N}}
\def\nut{\widetilde{\nu}}
\def\tfrac{\textstyle\frac}
\def\b{\beta}
\def\G{\Gamma}
\def\g{\gamma}
\def\zt{{\Bbb Z}_2}
\def\zth{{\bold Z}_2^\wedge}
\def\bs{{\bold s}}
\def\bx{{\bold x}}
\def\bof{{\bold f}}
\def\bq{{\bold Q}}
\def\be{{\bold e}}
\def\lline{\rule{.6in}{.6pt}}
\def\xb{{\overline x}}
\def\xbar{{\overline x}}
\def\ybar{{\overline y}}
\def\zbar{{\overline z}}
\def\ebar{{\overline \be}}
\def\nbar{{\overline n}}
\def\rbar{{\overline r}}
\def\Mbar{{\overline M}}
\def\et{{\widetilde e}}
\def\ni{\noindent}
\def\ms{\medskip}
\def\ehat{{\hat e}}
\def\xhat{{\widehat x}}
\def\nbar{{\overline{n}}}
\def\minp{\min\nolimits'}
\def\mul{\on{mul}}
\def\N{{\Bbb N}}
\def\Z{{\Bbb Z}}
\def\Q{{\Bbb Q}}
\def\R{{\Bbb R}}
\def\C{{\Bbb C}}
\def\notint{\cancel\cap}
\def\el{\ell}
\def\TC{\on{TC}}
\def\dstyle{\displaystyle}
\def\ds{\dstyle}
\def\Remark{\noindent{\it  Remark}}
\title
{On the cohomology classes of planar polygon spaces}
\author{Donald M. Davis}
\address{Department of Mathematics, Lehigh University\\Bethlehem, PA 18015, USA}
\email{dmd1@lehigh.edu}
\date{April 14, 2016}


\maketitle
\begin{abstract} We obtain an explicit formula for the Poincar\'e duality isomorphism $H^{n-3}(\Mbar(\ell);\zt)\to\zt$ for the space of isometry classes of $n$-gons with specified side lengths, if $\ell$ is monogenic in the sense of Hausmann-Rodriguez. This has potential application to topological complexity.
 \end{abstract}

\section{Main theorem}\label{intro}
If $\ell=(\ell_1,\ldots,\ell_n)$ is an $n$-tuple of positive real numbers, let $\Mbar(\ell)$ denote the space of isometry classes of oriented $n$-gons in the plane with the prescribed side lengths. In \cite{HK}, a complete description of $H^*(\Mbar(\ell);\zt)$ was given in terms of generators and a complicated set of relations.  In \cite{D}, explicit calculations were made in $H^*(\Mbar(\ell);\zt)$ for length vectors $\ell$ satisfying certain conditions, enabling us to prove that, for these $\ell$, the topological complexity of $\Mbar(\ell)$ satisfied
\begin{equation}\label{TC}2n-6\le\TC(\Mbar(\ell))\le 2n-5.\end{equation}
This is a result in topological robotics, as it specifies the number of motion planning rules required for a certain $n$-armed robot.(\cite{F}) However, our result only applied to a very restricted set of length vectors $\ell$.

The groups $H^k(\Mbar(\ell);\zt)$ are spanned by monomials $R^{k-r}V_{j_1}\cdots V_{j_r}$ for distinct positive subscripts $j\le n-1$. Here $R$ and $V_j$ are elements of $H^1(\Mbar(\ell);\zt)$. Since $\Mbar(\ell)$ is an $(n-3)$-manifold, there is a Poincar\'e duality isomorphism
\begin{equation}\label{PD}\phi:H^{n-3}(\Mbar(\ell);\zt)\to\zt.\end{equation}
For the cases considered in \cite{D}, we obtained an explicit formula for $\phi(R^{n-3-r}V_{j_1}\cdots V_{j_r})$. The contribution of this paper is to extend that formula to a broader class of length vectors. Note that it tells for each monomial whether it is 0 or the nonzero class, hence the title.

In order to describe these length vectors, we review the notion of genetic code introduced in \cite{H}. Since permuting the length vectors does not affect the homeomorphism type of $\Mbar(\ell)$, we may assume that $\ell_1\le\cdots\le\ell_n$.  A subset $S\subset\{1,\ldots,n\}$ is called {\it short} if $\sum_{i\in S}\ell_i<\sum_{i\not\in S}\ell_i$. We assume that $\ell$ is {\it generic}, which says that there are no subsets $S$ for which $\sum_{i\in S}\ell_i=\sum_{i\not\in S}\ell_i$. For generic $\ell$, a subset which is not short is called {\it long}. We define a partial order on sets  of integers by $\{s_1,\ldots,s_k\}\le T$ if there exist distinct $t_1,\ldots,t_k$ in $T$ with $s_i\le t_i$. The {\it genetic code} of $\ell$ is the set of maximal elements in the set of short subsets of $\ell$ which contain $n$. An element in the genetic code is called a {\it gene}.

One of the main theorems of \cite{D} was that, with three exceptions, (\ref{TC}) holds if $\ell$ has a single gene of size $4$. In order to prove this, we needed and obtained  the explicit formula for
$\phi(R^{n-3-r}V_{j_1}\cdots V_{j_r})$ for such length vectors.(\cite[Thm 4.1]{D}) In this paper, we extend this formula to all monogenic codes. We hope that this formula will enable us to study the cohomological implications for topological complexity of these spaces. The huge variety of genetic codes makes it seem unlikely that a formula such as ours might be extended to all genetic codes.

In \cite{D}, we introduced the term {\it gee} to refer to a gene with the $n$ omitted. This is sensible since, for all genes $G$, $n\in G$. Also, most of the formulas do not involve $n$. We say that a {\it subgee} is any set of positive integers which is less than or equal to a gee. Thus the subgees are all sets $S$ for which $S\cup\{ n\}$ is short. It is customary to write the elements of a gee in decreasing order. Only those $\{j_1,\ldots,j_r\}$  which are subgees can have $R^{k-r}V_{j_1}\cdots V_{j_r}\ne0$.(\cite[Cor 9.2]{HK})

Now we can state our theorem. For $B=(b_1,\ldots,b_k)$, let $|B|=\sum b_i$.
\begin{thm}\label{mainthm} Suppose $\ell$ has a single gee, $G=\{a_1+\cdots+a_k,\ldots,a_1+a_2,a_1\}$ with $a_i>0$. If $J$ is a set of distinct integers $\le a_1+\cdots+ a_k$, let $\theta(J)=(\theta_1,\ldots,\theta_k)$, where $\theta_i$ is the number of elements  $j\in J$ satisfying
$$a_1+\cdots+a_{i-1}<j\le a_1+\cdots+a_i.$$
Let $\cS_k$ denote the set of $k$-tuples of nonnegative integers such that, for all $i$, the sum of the last $i$ components of the $k$-tuple is $\le i$.
Then, for $\phi$ as in (\ref{PD}),
\begin{equation}\label{maineq}\phi(R^{n-3-r}V_{j_1}\cdots V_{j_r})=\sum_B\prod_{i=1}^k\tbinom{a_i+b_i-2}{b_i},\end{equation}
where $B$ ranges over all $(b_1,\ldots,b_k)$ for which $|B|=k-r$ and $B+\theta(\{j_1,\ldots,j_r\})\in\cS_k$.\end{thm}

Note that $J=\{j_1,\ldots,j_r\}$ is a subgee if and only if $\theta(J)\in\cS_k$. For example, if two $j$'s are greater than $a_1+\cdots+a_{k-1}$, then $\theta(J)$ has last component greater than 1, so $J$ is not in $\cS_k$, and $J$ cannot be $\le G$ since it has two elements greater than the second largest element of $G$. Thus (\ref{maineq}) is only relevant when $J$ is a subgee, but it yields 0 in other cases, anyway.

We illustrate Theorem \ref{mainthm} with several examples. First note that if $r=k$, then $\phi(R^{n-3-r}V_{j_1}\cdots V_{j_r})=1$. It seems that one cannot prove this simple and useful formula without deriving the entire formula for all terms.
Proving this case was a major motivation for this paper.

 In working with our formula,  it is useful to denote by $Y_T$ any term $R^{n-3-r}V_{j_1}\cdots V_{j_r}$
for which $\theta(\{j_1,\ldots,j_r\})=T$. The reader can verify that the case $k=3$ of Theorem \ref{mainthm} agrees with \cite[Thm 4.1]{D}, when the latter is expressed as in Example \ref{expl}, writing $a_i'$ for $a_i-1$. For example, in $\phi(Y_{0,1,0})$, $B=(2,0,0)$, $(1,1,0)$, and $(1,0,1)$ satisfy $B+(0,1,0)\in\cS_3$, but $B=(0,1,1)$ does not, since the sum of the last two entries of $(0,2,1)$ is greater than 2. (Our method of subscripting $Y$ here differs from that in \cite{D}.) This $\phi(Y_{0,1,0})$ refers to $\phi(R^{n-4}Y_j)$ for $a_1<j\le a_1+a_2$.

\begin{expl}\label{expl} If  $\ell$ has a single gee $\{a_1+a_2+a_3,a_1+a_2,a_1\}$, then
\begin{eqnarray*}
\phi(Y_T)&=&1\text{ if }|T|=3\text{ and }T\in\cS_3\\
\phi(Y_{0,2,0})=\phi(Y_{0,1,1})&=&a_1'\\
\phi(Y_{1,0,1})&=&a_1'+a_2'\\
\phi(Y_{2,0,0})=\phi(Y_{1,1,0})&=&a_1'+a_2'+a_3'\\
\phi(Y_{0,0,1})&=&\tbinom {a_1}2+a_1'a_2'\\
\phi(Y_{0,1,0})&=&\tbinom {a_1}2+a_1'a_2'+a_1'a_3'\\
\phi(Y_{1,0,0})&=&\tbinom {a_1}2+\tbinom {a_2}2+a_1'a_2'+a_1'a_3'+a_2'a_3'\\
\phi(Y_{0,0,0})&=&\tbinom {a_1}2(a_1'+a_2'+a_3')+\tbinom {a_2}2a_1'+a_1'a_2'a_3'.
\end{eqnarray*}
\end{expl}

\section{Proof}
In this section, we prove Theorem \ref{mainthm}.  As noted above, $H^{n-3}(\Mbar(\ell);\zt)$ is spanned by monomials $R^{n-3-r}V_{j_1}\cdots V_{j_r}$ for which $J=\{j_1,\ldots,j_r\}\le G$, i.e., $J$ is a subgee. Using \cite[Cor 9.2]{HK} as interpreted in \cite[Thm 2.1]{D}, a complete set of relations
is given by relations $\cR_I$ for each subgee $I$ except the empty set. This relation $\cR_I$ says
\begin{equation}\label{reln}\sum_{J\notint I}R^{n-3-|J|}\prod_{j\in J}V_j=0,\end{equation}
where the sum is taken over all subgees $J$ disjoint from $I$. To prove our theorem, it suffices to show that our proposed $\phi$ sends each relation $\cR_I$ to 0, since it will then be the unique nonzero homomorphism (\ref{PD}).

Similarly to \cite[(4.2)]{D}, the number of subgees $J$ disjoint from $I$ and satisfying $\theta(J)=C=(c_1,\ldots,c_k)$ is $\ds\prod_{i=1}^k\tbinom{a_i-m_i}{c_i}$, if $\theta(I)=(m_1,\ldots,m_k)$. Note that $m_i\le a_i$ since $m_i$ is the size of a subset of $a_i$ integers. Since our $\phi$ applied to a term in (\ref{reln}) is determined by $\theta(J)$, $\phi(\cR_I)$ becomes
\begin{equation}\label{formula}\sum_C\prod_{i=1}^k\tbinom{a_i-m_i}{c_i}\sum_B\prod_{i=1}^k\tbinom{a_i+b_i-2}{b_i},\end{equation}
where $|B|=k-|C|$ and $B+C\in\cS_k$. Letting $T=B+C$, this can be rewritten, with the first sum taken over all $T\in \cS_k$ satisfying $|T|=k$, as
\begin{eqnarray}&&\sum_T\sum_{B\le T}\prod_{i=1}^k\tbinom{a_i-m_i}{t_i-b_i}\tbinom{a_i+b_i-2}{b_i}\nonumber\\
&=&\sum_T\prod_{i=1}^k\sum_{b_i}\tbinom{a_i-m_i}{t_i-b_i}\tbinom{a_i+b_i-2}{b_i}\nonumber\\
&\equiv&\sum_T\prod_{i=1}^k\sum_{b_i}\tbinom{a_i-m_i}{t_i-b_i}\tbinom{1-a_i}{b_i}\nonumber\\
&=&\sum_T\prod_{i=1}^k\tbinom{1-m_i}{t_i}.\label{eq4}\end{eqnarray}
Here we use that $\binom{a+b-2}b=\pm\binom{1-a}b$.

The sum (\ref{eq4}) can be considered as a sum $\Sigma_1$ over {\it all} $k$-tuples $T$ of nonnegative integers summing to $k$ minus the sum $\Sigma_2$ over those which are not in $\cS_k$.
The sum $\Sigma_1$ equals $\binom{k-\sum m_i}k$, which is 0 unless $\sum m_i=0$, but that case has been excluded. ($I\ne\emptyset$.)

If $t_1,\ldots,t_j$ are nonnegative integers for which $t_1+\cdots+t_j<j$, but  $t_1+\cdots+t_i\ge i$ for all $i<j$, let $U(t_1,\ldots,t_j)$ denote the set of $k$-tuples $T$ indexing the sum $\Sigma_2$ which begin with $(t_1,\ldots,t_j)$. Since $t_1+\cdots+t_j<j$ is equivalent to saying that the sum of the last $k-j$ components is not $\le k-j$, these sets $U(t_1,\ldots,t_j)$ partition the set of $T$'s which occur in $\Sigma_2$. We show that the sum $\Sigma_2$ restricted to any such set $U(t_1,\ldots,t_j)$ is 0, which will complete the proof. We have
\begin{eqnarray*}&&\sum_{T\in U(t_1,\ldots,t_j)}\prod_{i=1}^k\tbinom{1-m_i}{t_i}\\
&=&\prod_{i=1}^j\tbinom{1-m_i}{t_i}\cdot\sum_{|T'|=k-t_1-\cdots-t_j}\prod_{i=j+1}^k\tbinom{1-m_i}{t_i}\\
&=&\prod_{i=1}^j\tbinom{1-m_i}{t_i}\cdot\tbinom{k-j-m_{j+1}-\cdots-m_k}{k-t_1-\cdots-t_j}.\end{eqnarray*}
In the second line, $T'=(t_{j+1},\ldots,t_k)$.
Since $(m_1,\ldots,m_k)$ arises from a subgee, $m_{j+1}+\cdots+m_k\le k-j$. But $k-t_1-\cdots-t_j>k-j$. Thus the final binomial coefficient consists of a nonnegative integer atop a larger integer, and hence is 0.

 \def\line{\rule{.6in}{.6pt}}


\begin{thebibliography}{99}
\bibitem{D} D.M.Davis, {\em Topological complexity of planar polygon spaces with small genetic code}, on {\tt arXiv}.
%{\em Topological complexity of some planar polygon spaces}, on arXiv.
%\bibitem{maple} \line, {\tt http://www.lehigh.edu/$\sim$dmd1/mapleresults.pdf}.
%\bibitem{27} \line, {\tt http://www.lehigh.edu/$\sim$dmd1/27cases.pdf}.
%\bibitem{Da} \line, {\em Immersions of real projective spaces}, Proceedings of Lefschetz Conference, Contemporary Math, Amer Math Soc. {\bf 58} (1987) 31--42.
\bibitem{F} M.Farber, {\em Invitation to topological robotics}, European Math Society (2008).
%\bibitem{F2} \line, {\em Topological complexity of motion planning}, Discrete Comput.~Geom. {\bf 29} (2003) 211--221.
%\bibitem{FTY} M.Farber, S.Tabachnikov, and S.Yuzvinsky, {\em Topological robotics: motion planning in projective spaces}, Intl Math  Research Notices {\bf 34} (2003) 1853--1870.
%\bibitem{Hb} J.-C.Hausmann, {\em Sur la topologie des bras articul\'es}, Lecture Notes in Mathematics, Springer {\bf 1474} (1989) 146--159.
\bibitem{HK} J.-C.Hausmann and A.Knutson, {\em The cohomology rings of polygon spaces}, Ann Inst Fourier (Grenoble) {\bf 48} (1998) 281--321.
\bibitem{H} J.-C.Hausmann and E.Rodriguez, {\em The space of clouds in Euclidean space}, Experiment Math {\bf 13} (2004) 31--47.
%\bibitem{web} \line, {\tt http://www.unige.ch/math/folks/hausmann/polygones}
%\bibitem{KK} Y.Kamiyama and K.Kimoto, {\em The height of a class in the cohomology ring of polygon spaces}, Int Jour of Math and Math Sci (2013) 7 pages.
%\bibitem{Kl} A.Klyachko, {\em Spatial polygons and stable configurations of points in the projective line}, Algebraic Geometry and its Applications (Yaroslav, 1992), Aspects Math, Vieweg, Braunschweig (1994) 67--84.
\end{thebibliography}
\end{document}